\documentclass[12pt]{amsart}
\usepackage{geometry}
\usepackage{amsmath,amsfonts,amssymb, tikz, amsthm}
\usepackage{hyperref}

\usepackage[utf8]{inputenc}
\usepackage{graphicx}
\usepackage[english]{babel}
\usepackage[tableposition=top]{caption} 
\usepackage{amsmath}
\usepackage{mathrsfs}
\usepackage{amssymb}
\usepackage{amsfonts}
\usepackage{mathtools}
\usepackage{tikz-cd}
\usepackage[utf8]{inputenc}

\newtheorem{theorem}{Theorem}[section]
\newtheorem{corollary}[theorem]{Corollary}
\newtheorem{lemma}[theorem]{Lemma}
\newtheorem{proposition}[theorem]{Proposition}
\newtheorem{definition}[theorem]{Definition}

\theoremstyle{remark}

\newtheorem*{rem}{Remark}

\theoremstyle{definition}
\newtheorem*{note}{Notation}
\newtheorem*{conv}{Convention}

\numberwithin{equation}{section}

\newcommand{\Z}{\mathbb{Z}}

\newcommand{\Q}{\mathbb{Q}}

\newcommand{\R}{\mathbb{R}}

\newcommand{\Qtr}{\Q^{\text{tr}}}
\newcommand{\N}{\text{N}}

\newcommand{\co}{\mathcal{O}}

\makeatletter
\makeatletter
\newcommand*{\house}[1]{%
  \mathord{%
    \mathpalette\@house{#1}%
  }%
}
\newcommand*{\@house}[2]{%
  \dimen@=\fontdimen8 %
      \ifx#1\scriptscriptstyle\scriptscriptfont
      \else\ifx#1\scriptstyle\scriptfont
      \else\textfont\fi\fi
      3 %
  \sbox0{%
    $#1%
      \vrule width\dimen@\relax
      \overline{%
        \kern2\dimen@
        \begingroup 
          #2%
        \endgroup
        \kern2\dimen@
      }%
      \vrule width\dimen@\relax
      \mathsurround=1.5\dimen@ 
    $%
  }%
  \ht0=\dimexpr\ht0-\dimen@\relax
  \dp0=\dimexpr\dp0+2\dimen@\relax
  \vbox{%
    \kern\dimen@ 
    \copy0 %
  }%
}

\begin{document}

\title{Northcott property and universality of higher degree forms}

\author{Om Prakash}
\address{Charles University, Faculty of Mathematics and Physics, Department of Algebra, Soko\-lo\-vsk\'{a} 83, 186 75 Praha 8, Czech Republic}
\email{prakash@karlin.mff.cuni.cz}

\keywords{Northcott property, $m$-ic lattices, universal higher degree forms, totally real fields, orders}
\subjclass[2020] {11E10, 11E76, 11G50, 11H06, 11P05, 11R04, 11R20, 11R80}

\thanks{The author was supported by {Czech Science Foundation} grant 21-00420M, {Charles University} programmes PRIMUS/24/SCI/010 and UNCE/24/SCI/022.}

\begin{abstract} 
Let $K$ be a totally real number field, $d$ a positive integer, and $Q$ a higher degree form over $K$. We prove that there are at most finitely many totally real extensions $L/K$ of degree $d$ such that $Q$ over $L$ is universal. Further, we show that there are no universal forms over totally real infinite extensions of $\Q$ having the Northcott property. 
\end{abstract}

\maketitle

\section{Introduction}
A \emph{universal quadratic form} over rational integers is a positive definite quadratic form that represents all positive rational integers. The $290$-Theorem of Bhargava-Hanke \cite{bh} fully characterizes universal quadratic forms over rational integers as exactly the forms that represent $1,2,\dots,290$. The situation becomes even more intriguing when we extend the study of universal forms over rational integers to the ring of integers $\co_K$ in a \emph{totally real} number field $K$. A quadratic form over a totally real number field $K$ is \emph{universal} if it is totally positive definite and represents all $\co_K^+$, i.e., all totally positive integers in $\co_K$. 

In this context, the estimates on small elements in $K$ have been quite effective, for example, the \emph{indecomposables}, i.e., those $\alpha\in\co_K^+$ that cannot be written as $\alpha=\beta+\gamma$ for some $\beta,\gamma\in\co_K^+$, have been proved very useful. Blomer-Kala \cite{bk1, bk2, k} showed how effectively the indecomposables can be used in estimating the minimal rank of universal quadratic forms over real quadratic fields. Another pertinent example is elements in $\co_K$ of small norm or trace. For instance, Kala \cite{k2} using elements of small norm extended the result of \cite{bk1}. By characterizing all the indecomposables and the elements of small trace in certain cubic fields, Yatsyna \cite{y} and Kala-Tinkov\'{a} \cite{kt} obtained an estimate on the minimal rank of universal quadratic forms over these cubic fields. For more applications of these tools, we direct readers to the survey \cite{survey} and the references therein. One of the reasons behind the importance of these small elements is that they are hard to be non-trivially represented by a quadratic form.

The theme of this article is to study the representation of algebraic integers by an \emph{$m$-ic form}, i.e., a homogeneous polynomial of degree $m\geq 2$, over totally real extensions of $\Q$. The main motivation behind this comes from the well known Waring's problem, which asks whether there exists a constant $g(m)$ such that every positive integer can be written as a sum of at most $g(m)$ $m$th powers of non-negative integers. In $1909$, Hilbert established the existence of $g(m)$ for each $m\geq 1$, moreover, due to the efforts of many mathematicians, e.g., see \cite{kw, vw}, we now almost know the value of $g(m)$, which is $g=2^m+\lfloor (3/2)^m\rfloor -2$ except for finitely many $m$. Later, Siegel \cite{si1, si2, siegel} studied Waring's problem over algebraic number fields. Over a number field $K$, Waring's problem can be stated as ``Does there exist a constant $G(K,m)$ such that every totally positive integer in $\sum\co_K^m$ (subring generated by $m$th powers of integers in $\co_K$) of sufficiently large norm can be written as a sum of at most $G(K,m)$ $m$th powers totally positive integers in $K$?" It is both intuitive and elegant to study integers represented by a higher degree form rather than just by sum of $m$th powers. 

In this paper we are interested in universal $m$-ic forms, i.e., totally positive definite $m$-ic forms that represent all totally positive integers in $\co_K$. Recently, Kala and the author \cite[Theorem $1.2$]{kp} showed that statements like $290$-Theorem would not hold for $m$-ic forms, i.e., universal $m$-ic forms cannot be characterized by a finite subset of totally positive integers. This is one illustration of the fact that the study of universal $m$-ic forms is more complicated than the quadratic case.

To understand the difficulty in the case of $m$-ic forms, note that, to a quadratic form we can associate the Gram matrix. This helps in understanding various properties, e.g., positive definiteness. Whereas, for $m$-ic forms we don't have such a matrix. There are several other foundational concepts that these $m$-ic forms lack. 

On the other hand, several results from the algebraic theory of quadratic forms have been extended to the $m$-ic setting. Beginning with Harrison's \cite{harrison} foundational work on Witt theory for $m$-ic forms over field of sufficiently large characteristic. Following this, various other results were also extended by several mathematicians (e.g., see \cite{ab1,ab2,ab3,ab4,ab5,ab6,pumpl}). 

In this article, we first implement the language of $m$-ic lattices over totally real extension of $\Q$, which is a slight generalization of $m$-ic forms and gives a coordinate free approach. 
One may encounter that not all totally positive integers can written as a sum of $m$th powers. Thus, it makes sense to study the representation of multiples of $\co_K^+$. Our first main result of this paper is the following: 
\begin{theorem}
    Let $m\in\Z_{\geq2}$ be even, $K$ be a totally real number field, $(\Lambda, Q)$ a totally positive definite $m$-ic $\co_K$-lattice, and $\alpha\in\co_K^+.$ There exist at most finitely many totally real number fields $L/K$ of degree $d=[L : \Q]$ such that $(\Lambda\otimes_{\co_K}\co_L, Q_L)$ represents all $\alpha\co_L^+.$
\end{theorem}

For quadratic forms, the above result was proved by Kala-Yatsyna \cite[Theorem $7$]{ky}. The key ingredient that led us to this is our Proposition \ref{p1}. 

We say that a quadratic form is \emph{isotropic} if it represents $0$ non-trivially. Springer's theorem for odd degree extensions states that, given a quadratic form $Q$ over $K$ and an extension $L/K$ of odd degree, then $Q$ is isotropic over $K$ if and only if it is isotropic over $L$. Johnson-Pumpl\"{u}n \cite{pumpl} discussed the possibility of extending Springer's theorem to higher degree forms. They verified the truth of Springer's theorem in case when the form permits composition or Jordan composition. In general, this is widely open. Recently, Daans-Kala-Kr\'{a}sensk\'{y}-Yatsyna \cite[Question]{dkky} asked the Springer's theorem in the context of integral representation, i.e., if an integer $\alpha\in\co_K$ is represented by the $Q$ over $L$, does it imply that $\alpha$ is also represented by $Q$ over $K$? In \cite[Theorem $3.1$]{dkky}, they showed that this may fail. By presenting another application of our Proposition \ref{p1}, we extend their result to the setting of $m$-ic forms in Theorem \ref{fail.ST}. 

\medskip 

In Section \ref{s6}, we turn our attention to a significantly different scenario involving $m$-ic forms over totally real infinite extensions of $\Q$. 
In a fairly recent work, Daans-Kala-Man \cite{dkm} pioneered the study of quadratic forms over totally real infinite extensions of $\Q$. They proposed two possible ways to study the universality of quadratic forms over these extensions. One method uses the Northcott property, while the other relies on the number of square classes of totally positive units. We adopt the first method to study the universality of $m$-ic forms over infinite extensions, and we prove the following result: 

\begin{theorem}\label{main2}
    Let $m\in\Z_{\geq 2}$ be even, $K$ be a totally real infinite extension of $\Q$, $\delta\in\co_K^+$, and $\co$ be an order in $K$. If there exists an $m$-ic $\co$-lattice that represents all $\delta\co^+$, then $\co^+$ does not have the Northcott property. 
\end{theorem}

Note that when $K$ is a totally real finite extension of $\Q$, Kala and the author \cite[Theorem $5.4$]{kp}, using Waring's problem and some estimates from geometry of numbers, have established that universal forms always exist.

Theorem \ref{main2} in fact gives a bunch of examples (see Corollary \ref{c2}) of totally real infinite extensions of $\Q$ that have no universal $m$-ic lattices. 

An extremely interesting question in Diophantine geometry is -- \emph{which infinite extension has the Northcott property?} Concerning this, Theorem \ref{main2} can also be viewed as a criterion for a totally real infinite extension of $\Q$ to have the Northcott property.  

\section*{Acknowledgements}
The author wishes to express his thanks to V\'{i}t\v{e}zslav Kala for many helpful comments and suggestions throughout the preparation of this article. The author would also like to thank Nicolas Daans and Pavlo Yatsyna for answering various questions.  

\section{Preliminaries}\label{s2}

\subsection{Finite extensions of $\Q$} Let $K$ be a \emph{totally real} number field of degree $d=[K:\Q]$, i.e., all its embeddings $\sigma_1, \sigma_2,\dots, \sigma_d:K\rightarrow \mathbb{C}$ have image in $\mathbb{R}$, and $\co_K$ its ring of algebraic integers. An integral basis for $K$ is a set of $\Z$-linearly independent elements $\omega_1, \omega_2, \dots, \omega_d \in \mathcal{O}_K$ such that $\mathcal{O}_K = {\mathbb{Z}}\omega_1 + \cdots + {\mathbb{Z}}\omega_d$. 

We say that an element $\alpha\in K$ is \emph{totally positive} if $\sigma_i(\alpha)>0$ for all $1\leq i\leq d$, and write $\alpha\succ\beta$ if $\alpha-\beta$ is totally positive. For a subset $S\subseteq K$, $S^+$ denotes the set of all totally positive elements of $S.$ The \emph{trace} of $\alpha\in K$ is $\text{Tr}_{K/\Q}(\alpha)=\sum_{i=1}^d \sigma_i(\alpha)$, and the \emph{norm} of $\alpha$ is $\N_{K/\Q}(\alpha)=\prod_{i=1}^d\sigma_i(\alpha)$.

\subsection{Infinite extension of $\Q$} Since in the second part of the paper we will focus on studying forms over infinite extensions of $\Q$, let us collect some background information on infinite extensions. 

Let $K/\Q$ be a field extension, an element $\alpha\in K$ is said to be \emph{algebraic} over $\Q$ if $\alpha$ is a root of some polynomial with coefficients in $\Q$; and $K/\Q$ is said to be \emph{algebraic} if every element of $K$ is algebraic over $\Q$. A typical example of such an extension would be $\Q(\sqrt{n}\mid n\in \Z_{>1}).$ 

In this paper, we always assume that $K/\Q$ is an algebraic extension. Then for every $\alpha\in K$, we have that $\Q(\alpha)$ is a number field.

An element $\alpha\in K$ is \emph{totally real} if $\alpha$ is totally real in the number field $\Q(\alpha)$, and $\alpha$ is \emph{totally positive} if $\alpha$ is totally positive in $\Q(\alpha)$. We say that $K$ is \emph{totally real extension of $\Q$} if all elements of $K$ are totally real. Observe that $K$ is a totally real extension of $\Q$ if and only if each subfield $\Q\subseteq E \subseteq K$ of finite degree $[E:\Q]$ is totally real. 
From now on, we will only work with totally real extensions.

Let $\Qtr$ denote the \emph{maximal} totally real extension of $\Q$, i.e., every totally real extension of $\Q$ is contained in $\Qtr$. For $\alpha, \beta\in\Qtr$, we write $\alpha\succ\beta$ if $\alpha\succ\beta$ in $\Q(\alpha, \beta)$. $\co_K$ will denote the ring of algebraic integers in $K$; and for any subset $S\subseteq K$, $S^+$ denotes the set of totally positive elements of $K$.

In this paper, our height function will be the house $\house{ \cdot } :K\rightarrow \R_{\geq 0}$ defined as follows: for each $\alpha\in K$ the house of $\alpha$ is $\house{\alpha}=\max_{\sigma}(|\sigma(\alpha)|)$, where the maximum is taken over all embeddings $\sigma: \Q(\alpha)\rightarrow \R$. So, for example, when $K$ is a number field of degree $d$, then the \emph{house} of $\alpha$ is $\house{\alpha}=\max_{1\leq i\leq d} |\sigma_i(\alpha)|$, where $\sigma_is$ are embeddings of $K$.

For any non-zero $\alpha\in\co_K$ we have $\house{\alpha}\geq 1$, since $1\leq |\N_{\Q(\alpha)/\Q}(\alpha)|$. We recall some standard properties of house, which we will be using throughout.

\begin{itemize}
         \item $\house{\alpha+\beta}\leq \house{\alpha}+\house{\beta}$ for all algebraic numbers $\alpha, \beta$.
        \item $\house{\alpha\beta}=\house{\alpha}\house{\beta}$ for all algebraic numbers $\alpha, \beta$.
        \item (Northcott property \cite[Theorem $1$]{north}) Let $c_1,c_2 >0$ be real numbers. There exist only finitely many algebraic integers $\alpha$ with $\house{\alpha}<c_1$ and $\deg \alpha <c_2$.
\end{itemize}

For $v=(\gamma_1,\gamma_2,\dots,\gamma_n)\in K^n$, we define the house of $v$ as $\house{v}=\max_{1\leq i\leq n}\house{\gamma_i}.$

The most commonly used height function is the Weil height. Although we only require it for Corollary \ref{c2}, we would like to mention it here. 

For $\alpha\in K$ the \emph{absolute logarithmic height} or \emph{logarithmic Weil height} of $\alpha$ is defined as $$h(\alpha)=\frac{1}{[\Q(\alpha):\Q]}\sum_{v\in M_K} [\Q(\alpha)_v:\Q_v] \log \max(1,|\alpha|_v),$$ where $M_K$ is the set of all valuations of $\Q(\alpha)$, $\Q(\alpha)_v$ and $\Q_v$ are completions of $\Q(\alpha)$ and $\Q$ with respect to $v$, and $|\alpha|_v$ is the normalized valuation on $\alpha$. In particular, when $\alpha$ is an algebraic integer in some totally real field $K$, we have $$h(\alpha)=\frac{1}{[\Q(\alpha):\Q]}\sum_\sigma \log \max(1,|\sigma(\alpha)|),$$ where the sum runs over all the embeddings $\sigma$ of $\Q(\alpha)$. The Northcott property with respect to Weil height then asserts that, for given real numbers $c_1,c_2>0$, there exist only finitely many $\alpha\in\co_K$ with $h(\alpha)<c_1$ and $\deg \alpha <c_2$.

If $0\neq \alpha\in\co_K$, then from \cite[Lemma $7$]{vv16} one can obtain that $h(\alpha)\leq \log \house{\alpha}$. Thus if $\co_K^+$ has the Northcott property with respect to Weil height, then it has Northcott property with respect to the house.

\begin{conv}
    Throughout the article, whenever we refer to the Northcott property, we always mean the Northcott property with respect to the house. 
\end{conv}

\section{Introduction to $m$-ic lattices}\label{s3}
Throughout this section, let $m$ be a positive integer $\geq 2$. At the beginning of this section, $R$ will be an integral domain of characteristic $0$. While most of the notions we define will be applicable in characteristics $p>m$, where $p$ is a rational prime, in the rest of the paper we will work only in characteristic $0$. Therefore, to avoid any confusion, we will define everything over characteristic $0$.

Let $F$ be the field of fractions of $R$ denoted by $F=\text{Frac }R$ and $V$ be a finite dimensional $F$-vector space. An $F$-\emph{multilinear form} is a map 
$$M : \underbrace{V\times V\times \cdots \times V}_{m-\text{times}} \rightarrow F$$
which is $F$-linear in each variable.

\begin{definition}\label{def: mic-space}
     Let $V$ be a finite dimensional $F$-vector space. An \emph{$m$-ic form} over $F$ is a map $Q : V \rightarrow F$ such that     
    \begin{enumerate}
        \item $Q(av)=a^m Q(v)$ for all $v\in V$, $a\in F,$ and
        \item the map $$M_Q : \underbrace{V\times V\times \cdots \times V}_{m-\text{times}} \rightarrow F$$ $$(v_1,\dots,v_m)\mapsto \frac{1}{m!}\sum_{I\subseteq \{1,2,\dots m\}} (-1)^{m-|I|} Q\left(\sum_{i\in I} v_i\right),$$ (with the convention that when $I=\emptyset$ then the sum is $0$) is an $F$-multilinear form.
    \end{enumerate}
    We call the pair $(V, Q)$ an $m$-ic space over $F$, and we define the dimension of $(V,Q)=\dim V$.
\end{definition}
Note that by definition, $M_Q$ is symmetric in the sense that $M_Q(v_1, \dots, v_m)=M_Q(\tau(v_1), \dots, \tau(v_m))$ for every permutation $\tau$. We can recapture $Q$ from $M_Q$ by setting $M_Q(v,\dots,v)=Q(v)$, thanks to the identity $m!=\sum_{r=1}^m (-1)^{m-r} \binom{m}{r}r^m$.
\medskip

Given $m$-ic spaces $(V_1,Q_1)$ and$(V_2,Q_2)$ over $F$, an isomorphism of $F$-vector space $\psi:V_1\rightarrow V_2$ is called an \emph{isometry between $(V_1,Q_1)$ and $(V_2,Q_2)$} if for all $v\in V_1$ we have $Q_1(v)=Q_2(\psi(v))$.

The following fact is well known (cf. \cite{harrison}, \cite{huang}). However, for the sake of the reader, let us prove it here. 
\begin{proposition}\label{1-1 space}
    Let $m\in\Z_{\geq 2},n\in\Z_{>0}$, and let $f\in F[x_1,x_2,\dots,x_n]$ be a homogeneous polynomial of degree $m$. The map 
    \begin{equation*}
        Q_f: F^n\rightarrow F : (a_1,\dots,a_n)\mapsto f(a_1,\dots,a_n)
    \end{equation*}
    is an $m$-ic form.
    
    Conversely, given an $m$-ic space $(V,Q)$ over $F$ of dimension $n$, there exists a homogeneous polynomial $f\in F[x_1,x_2,\dots,x_n]$ of degree $m$ such that $(V,Q)$ is isometric to $(F^n,Q_f)$.
\end{proposition}

\begin{proof}
    Let $f(x_1,x_2,\dots,x_n)=\sum_{1\leq i_1,\cdots,i_m\leq n}a_{i_1\cdots i_m}x_{i_1}\cdots x_{i_m}\in F[x_1,x_2,\dots,x_n]$ be arbitrary homogeneous polynomial of degree $m$.
    To show $Q_f$ is an $m$-ic form, we need to verify that the map $Q_f$ satisfies the conditions stated in Definition \ref{def: mic-space}. Let $e_1,\dots,e_n$ be an $F$-basis for $F^n$, then every $v\in F^n$ may be represented by the coordinates $(x_1,\dots,x_n)\in F^n$. 

    Since $f$ is homogeneous, observe that for all $v=(\alpha_1, \dots,\alpha_n)\in F^n$, and $a\in F$ we have that 
    $$Q_f(av)=Q_f(a\alpha_1,\dots,a\alpha_n)=f(a\alpha_1,\dots,a\alpha_n)=a^nf(\alpha_1,\dots,\alpha_n)=a^nQ_f(\alpha_1,\dots,\alpha_n)=a^nQ_f(v).$$

    Now, it remains to show that the map $$M_{Q_f} :  \underbrace{F^n\times \cdots \times F^n}_{m-\text{times}} \rightarrow F$$
    $$(v_1,\dots,v_m)\mapsto \frac{1}{m!}\sum_{I\subseteq \{1,2,\dots m\}} (-1)^{m-|I|} Q_f\left(\sum_{i\in I} v_i\right)$$
    is a $F$-multilinear form. Let us write $v_i=(x_1^{(i)},\dots,x_n^{(i)})$. Since $M_{Q_f}$ is symmetric, it is enough to verify linearity in the first variable. 

    Let $\alpha\in F$ and $u_1=(y_1^{(1)},\dots,y_n^{(1)})\in F^n$ be arbitrary. In the following, $I_1$ is a subset of $\{1,2,\dots,m\}$ containing $1$, and $I_{-1}$ is a subset of $\{1,2,\dots,m\}$ not containing $1$. We then have
    \begin{equation*}
        \begin{split}
             m! M_{Q_f}(\alpha v_1+u_1,v_2,\dots,v_m) & =
             \sum_{I_1\subseteq \{1,2,\dots m\}} (-1)^{m-|I_1|} Q_f\left(\left(\alpha v_1+u_1\right)+\sum_{i\in I_1\setminus \{1\}} v_i\right) \\
             & + \sum_{I_{-1}\subseteq \{2,\dots m\}} (-1)^{m-|I_{-1}|} Q_f\left(\sum_{i\in I_{-1}\setminus \{1\}} v_i\right) \\
             & = \sum_{1\leq x_1,\cdots,i_m\leq n} a_{i_1\cdots i_m} \sum_{I_1\subseteq \{1,2,\dots,m\}}\left(-1\right)^{m-|I_1|} \prod_{j=1}^m\left(\alpha x_{i_j}^{(1)}+\sum_{\ell\in I_1\setminus \{1\}} x_{i_j}^{(\ell)}\right) \\
            & + \sum_{1\leq x_1,\cdots,i_m\leq n} a_{i_1\cdots i_m} \sum_{I_{-1}\subseteq \{2,\dots,m\}}\left(-1\right)^{m-|I_{-1}|} \prod_{j=1}^m\left(\sum_{\ell\in I_{-1}} x_{i_j}^{(\ell)}\right) \\
             & = \alpha \sum_{1\leq x_1,\cdots,i_m\leq n} a_{i_1\cdots i_m} \sum_{I\subseteq \{1,2,\dots,m\}}\left(-1\right)^{m-|I|} \prod_{j=1}^m\left(\sum_{\ell\in I} x_{i_j}^{(\ell)}\right)\\
             & + \sum_{1\leq x_1,\cdots,i_m\leq n} a_{i_1\cdots i_m} \sum_{I_1\subseteq \{1,2,\dots,m\}}\left(-1\right)^{m-|I_1|} \prod_{j=1}^m\left(y_{i_j}^{(1)}+\sum_{\ell\in I_1\setminus \{1\}} x_{i_j}^{(\ell)}\right) \\
            & + \sum_{1\leq x_1,\cdots,i_m\leq n} a_{i_1\cdots i_m} \sum_{I_{-1}\subseteq \{2,\dots,m\}}\left(-1\right)^{m-|I_{-1}|} \prod_{j=1}^m\left(\sum_{\ell\in I_{-1}} x_{i_j}^{(\ell)}\right) \\
            & = m!\alpha M_{Q_f}(v_1,\dots,v_m)+m! M_{Q_f}(u_1,v_2,\dots,v_m). 
        \end{split}
    \end{equation*}

    Conversely, let $(V,Q)$ be an $m$-ic space of dimension $n$, and $M_Q$ be the associated symmetric multilinear form. Let $u_1,u_2,\dots,u_n$ be an $F$-basis for $V$. Consider the map $$\psi : V\rightarrow F^n:v=x_1u_1+x_2u_2+\cdots+x_nu_n\mapsto (x_1,\dots,x_n).$$ It is straightforward to verify that $\psi$ is an isomorphism of $F$-vector spaces. Now, our aim is to show that there exists a homogeneous polynomial $f\in F[x_1,x_2,\dots,x_n]$ of degree $m$ such that $Q(v)=Q_f(\psi(v))$ for all $v\in V$, i.e., $Q(v)=Q_f(x_1,\dots,x_n)=f(x_1,\dots,x_n)$. 
    \medskip

    For each $1\leq j\leq m$, write $v_j=x_{1,j}u_1+x_{2,j}u_2+\cdots+x_{n,j}u_n$, where $x_{i,j}\in F$ for all $i$. Since $M_Q$ is symmetric multilinear form, we can write $M_Q$ as follows: 
    \begin{equation}\label{all multi}
       \begin{split}
            M_Q(v_1,\dots,v_m) & =\sum_{(i_1,\dots,i_m)\in\{1,\dots,n\}^m} x_{i_1,1}\cdots x_{i_m,m} M_Q(u_{i_1},\dots,u_{i_m})\\
            & = \underset{1\leq i,j\leq n}{\underset{I_i\cap I_j=\emptyset \text{ for }i\neq j}{\underset{I_1\cup I_2\cup\cdots\cup I_n = \{1,2,\dots,m\}}{\sum}}} \left( \prod_{i_1\in I_1} x_{1,i_1}\right)\cdots \left( \prod_{i_n\in I_n} x_{n,i_n}\right) M_Q(\underbrace{u_1,\dots,u_1}_{|I_1|}, \dots, \underbrace{u_n,\dots,u_n}_{|I_n|}).
       \end{split}
    \end{equation}
    Now, when we evaluate $Q(v)=M_Q(v,\dots,v)$, we obtain the same product $x_1^{k_1}\cdots x_n^{k_n}$ multiple times. This occurs as many times as there are ways to choose $n$ disjoint sets $I_j$, each with cardinality $k_i$, where $k_1+\cdots+k_n=m$, and this number is precisely $\binom{m}{k_1,\dots,k_n}=m!/k_1!\cdots k_n!$. Finally, we obtain: 
    \begin{equation*}
        Q(v)=\underset{0\leq k_i, 1\leq i\leq n}{\underset{k_1+\cdots+k_n=m}{\sum}} \binom{m}{k_1,\dots,k_n} x_1^{k_1}\cdots x_n^{k_n} M_Q(\underbrace{u_1,\dots u_1}_{k_1},\dots,\underbrace{u_n,\dots,u_n}_{k_n})=f(x_1,\dots,x_n), 
    \end{equation*}
    a homogeneous degree $m$ polynomial over $F$ in variables $x_1,\dots,x_n$. Hence $(V,Q)$ is isometric to $(F^n, Q_f)$.
\end{proof}

\begin{rem}
    In the converse part of the above proposition, we did not use the definition of $M_Q$ provided in Definition \ref{def: mic-space}. Therefore, there is a one-to-one correspondence between homogeneous polynomials of degree $m$ and multilinear forms. Also, note that any multilinear form $M:V\times\cdots\times V\rightarrow F$ is given by \ref{all multi}.
\end{rem}

An $R$-lattice is a finitely generated $R$-submodule $\Lambda\subseteq V$. The rank of $\Lambda$ is the $F$-dimension of $F$-span of $\Lambda$. A set of vectors $v_1,\dots,v_n\in \Lambda$ is called a \emph{free $R$-basis} for $\Lambda$ if it is $R$-independent and $Rv_1+Rv_2+\cdots+Rv_n=\Lambda$. A lattice that has a basis is called \emph{free}. Any two bases of a free lattice contain the same number of elements, and this number is the rank of $\Lambda$. It should be noted that not all lattices are free, for example, an $\co_K$-lattice $\Lambda$ is not free, when $K$ is a number field with class number $\neq 1$ \cite[$81:5$]{omeara}.
\begin{definition}
    Let $R$ be an integral domain with $\text{Frac }R=F$, $(V,Q)$ an $m$-ic space over $F$, $\Lambda\subseteq V$ an $R$-lattice. If we endow $\Lambda$ with the restriction map $Q_{|\Lambda} : \Lambda \rightarrow F$, then the pair $(\Lambda, Q_{|\Lambda})$ is called an $m$-ic $R$-lattice.
\end{definition}

\begin{note}
    For the rest of the article, we denote $Q_{|\Lambda}$ simply as $Q$. 
\end{note} 

We say that two $m$-ic $R$-lattices $(\Lambda_1, Q_1)$ and $(\Lambda_2, Q_2)$ are isometric if there exists an isomorphism $\psi: \Lambda_1 \rightarrow \Lambda_2$ of $F$-vector spaces such that $Q_1(v)=Q_2(\psi(v))$ for all $v\in \Lambda_1$.

\begin{proposition}
    Let $m\in\Z_{\geq 2},n\in\Z_{>0}$, $R$ be an integral domain with $F= \text{Frac } R$, and let $f\in F[x_1,x_2,\dots,x_n]$ be a homogeneous polynomial of degree $m$. The map 
    \begin{equation*}
        Q_f: F^n\rightarrow F : (a_1,\dots,a_n)\mapsto f(a_1,\dots,a_n)
    \end{equation*}
    is an $m$-ic form on $F^n$.
    
    Conversely, given a free $m$-ic $R$-lattice $(\Lambda, Q)$ of rank $n$, there exists a homogeneous polynomial $f\in F[x_1,x_2,\dots,x_n]$ of degree $m$ such that $(\Lambda, Q)$ is isometric to $(F^n,Q_f)$.
\end{proposition}

\begin{proof}
    Analogous to the proof of Proposition \ref{1-1 space}. 
\end{proof}

Let $R\subset S$ be integral domains. The tensor product $V\otimes_R S$ of $R$-modules $V$ and $S$ is an $R$-module equipped with a bilinear map $\otimes : V \times S\rightarrow V\otimes_R S: (v,\alpha)\mapsto v\otimes \alpha$. The universal property (up to $R$-module isomorphism) of this tensor product is that for any $R$-module $T$ and morphism $f: V \times S\rightarrow T$, there exists a unique morphism $\varphi:V \otimes_R S\rightarrow T$ such that $f=\varphi\circ\otimes$.  
We define the tensor product $(\Lambda, Q)\otimes_R S=(\Lambda\otimes_R S, Q_S)$, where $\Lambda\otimes_R S$ is tensor product of $R$-modules $\Lambda$ and $S$ and $Q_S: \Lambda\otimes_R S\rightarrow S$ is given by $Q_S(v\otimes \alpha)=\alpha^mQ(v)$. Thus, every $m$-ic $R$-lattice $(\Lambda, Q)$ can be associated to an $m$-ic $S$-lattice, namely $(\Lambda, Q)\otimes_R S=(\Lambda\otimes_R S, Q_S)$. We call $(\Lambda\otimes_R S, Q_S)$ an \emph{scalar extension} of $(\Lambda, Q)$.
\medskip 

In this article, we are interested in the situation where $F$ is a totally real extension of $\Q$, and $R\subseteq \co_F$ a subring such that $\text{Frac }R=F$.

\medskip

Let $K$ be a totally real extension of $\Q$ (not necessarily of finite degree). An \emph{order} in $K$ is a subring $\co \subseteq \co_K$ such that $\text{Frac } \co=K$. Let $V$ be a finite dimensional $K$-vector space. An $\co$-lattice is a finitely generated $\co$-submodule $\Lambda\subseteq V$; and the rank of $\Lambda$ is the $K$-dimension of $K$-span of $\Lambda.$ Let $Q: V \rightarrow K$ be an $m$-ic form. If we endow $\Lambda$ with the restriction map $Q_{|\Lambda} : \Lambda \rightarrow K$, then the pair $(\Lambda, Q_{|\Lambda})$ is called an $m$-ic $\co$-lattice. We say that an $m$-ic $\co$-lattice is integral if $Q(\Lambda)\subseteq \co$. Further, an integral lattice is said to be classical if $M_Q(v_1,\dots,v_m)\in \co$ for all $v_i\in \Lambda$. An element $\alpha\in K$ is represented by $(\Lambda,Q)$ if there exists $v\in \Lambda$ such that $Q(v)=\alpha$.

In this paper, we will focus on \emph{totally positive definite} $m$-ic lattices.

\begin{definition} \label{pos def}
Let $K$ be a totally real extension of $\Q$ and $(V,Q)$ be an $m$-ic space over $K$. Then $(V,Q)$ is positive definite if $Q(v)>0$ for all $v\in V\otimes_{\Q} \R.$ Further, $(V,Q)$ is said to be totally positive definite if for all embeddings $\sigma : K\rightarrow \R$, $(V\otimes_{\Q} \R, \sigma(Q))$, where $\sigma(Q)=\sigma \circ Q: V\rightarrow K\rightarrow \R$ is positive definite. 
\end{definition}

\begin{rem}
    The same definition as in the case of \emph{quadratic lattice}, i.e., $(V,Q)$ is positive definite if $Q(v)>0$ for all $v\in V$, would not be equivalent with Definition \ref{pos def} in the $m$-ic setting. For example, consider the form $$Q_f: \Q^2 \rightarrow \Q: (x,y)\mapsto f(x,y)=(x^2-2y^2)^2,$$ satisfies the criterion $Q(v)>0$ for all $v\in \Q^2$ but is not positive definite because $f(\sqrt{2}, 1)=0$ and $(\sqrt{2},1)$ is a non-zero element in $\Q^2\otimes_{\Q}\R$.
\end{rem}

From the above definition, it follows immediately that if $(V,Q)$ is positive definite then $Q(v)>0$ for all $v\in V$. 
Additionally, note that for $u\otimes -1\in V\otimes_{\Q} \R$, $Q(u\otimes -1)=(-1)^m Q(u)$. Therefore, whenever $Q$ is totally positive definite, we must have $m$ be even.

An $m$-ic $\co$-lattice $(\Lambda, Q)$ is totally positive definite if $(V,Q)$ is totally positive definite. We say that $(\Lambda, Q)$ is \emph{universal} if it is totally positive definite and represents all totally positive integers in $\co$. 

\begin{conv}
    By an $m$-ic $\co$-lattice, we always mean an integral $m$-ic $\co$-lattice.
\end{conv}

Let $L/K$ be an extension of totally real fields. Suppose $\co\subseteq \co_K$ and $\co \subseteq \co'\subseteq \co_L$ be orders in $K$ and $L$, i.e., $\text{Frac }\co = K$ and $\text{Frac }\co'=L.$ Recall that every $m$-ic $\co$-lattice can be extended to an $m$-ic $\co'$-lattice, namely, the scalar extension $(\co' \otimes_{\co} \Lambda, Q_L)$.

Conversely, given a subring $\co \subseteq \co_L$ such that $\text{Frac } \co =L$ and an $m$-ic $\co$-lattice $(\Lambda, Q)$, we say that $\Lambda$ is \emph{defined over} $K$ if $\co_{|K}=\co \cap \co_K$ satisfies $\text{Frac }\co_{|K} =K$, and there is an $\co_{|K}$-lattice $(\Lambda_{\co_{|K}}, Q_K)$ such that $(\Lambda, Q)$ is its scalar extension to $\co$.
\medskip

The following proposition is a straightforward generalization of \cite[Proposition $2.1$]{dkm}. However, for the sake of completeness, we will redo the proof in the context of $m$-ic lattices. 
\begin{proposition}\label{definedover}
    Let $n\in\Z_{>0}, m\in\Z_{\geq 2}$, $L/\Q$ be an algebraic extension, $\co$ an order in $L$, and $(\Lambda, Q)$ an $m$-ic $\co$-lattice of rank $n$. Then $(\Lambda, Q)$ is defined over some number field $K\subseteq L.$
\end{proposition}

\begin{proof}
    Let $(\Lambda, Q)$ be an $m$-ic $\co$-lattice of rank $n$. Suppose that $v_1,v_2,\dots,v_k\in\Lambda\subseteq V$ is a set of generators for $\Lambda$, where $k\geq n$. Without loss of generality assume that $v_1,v_2,\dots,v_n$ forms an $L$-basis for $V$. Then, for each $n< i\leq k$ we can write $$v_i=\gamma_{i1}v_1+\gamma_{i2}v_2+\cdots+\gamma_{in}v_n, \text{ where } \gamma_{ij}\in L.$$ Moreover, using the fact that $\co$ is an order in $L$, we may write $\gamma_{ij}=a_{ij}/b_{ij}$ with $a_{ij}, b_{ij}\in\co$. Write $M_Q(v_{i_1},v_{i_2}, v_{i_m})=\beta_{i_1,i_2,\dots,i_m}\in L$, then we can again write $\beta_{i_1,i_2,\dots,i_m}=c_{i_1,i_2,\dots,i_m}/d_{i_1,i_2,\dots,i_m}$ with $c_{i_1,i_2,\dots,i_m}, d_{i_1,i_2,\dots,i_m}\in \co$. 

    Let $K=\Q(a_{ij}, b_{ij},c_{i_1,i_2,\dots,i_m}, d_{i_1,i_2,\dots,i_m})\subset L$ be the subfield generated by $a_{ij}, b_{ij},c_{i_1,i_2,\dots,i_m}, d_{i_1,i_2,\dots,i_m}\in\co$. Note that $K$ is a number field since there are only finitely many $a_{ij}, b_{ij},c_{i_1,i_2,\dots,i_n}, d_{i_1,i_2,\dots,i_n}$, and it follows that $\co_{|_K}=\co\cap\co_K$ satisfies $\text{Frac }\co_{|_K}=K$. 

    Let $\Lambda_{\co_{|K}}$ be the $\co_{|K}$-submodule generated by $v_1,v_2,\dots,v_n$, $V_K$ be the $K$-vector subspace of $V$ generated by $v_1,v_2,\dots,v_n$, and $Q_K=Q_{|K}: \Lambda_{\co_{|K}}\rightarrow K$. Note that $\dim_K V_K =n$. 

    Consider $\co_{|K}$-modules $\Lambda_{\co_{|K}}\otimes_{\co_{|K}} L$ and $V_K\otimes_K L$. We have the map 
    $$f_1: \Lambda_{\co_{|K}}\times L\rightarrow V_K\otimes_K L: (v,\alpha)\mapsto v\otimes \alpha,$$
    which is a bilinear map, and the tensor product map
    $$\otimes: \Lambda_{\co_{|K}}\times L\rightarrow \Lambda_{\co_{|K}}\otimes_{\co_{|K}} L: (v, \alpha)\mapsto v\otimes \alpha.$$

    By universal property of tensor product, there exists a unique map $$\varphi_1: \Lambda_{\co_{|K}}\otimes_{\co_{|K}} L \rightarrow V_K\otimes_K L$$ such that $\varphi_1\circ \otimes=f_1,$ which gives $\varphi_1(v\otimes\alpha)=v\otimes\alpha$. 

    Recall that $V_K$ is generated by $v_1,\dots,v_n$. Thus, we can define the following  map on generators:   
    $$f_2 : V_K\times L\rightarrow \Lambda_{\co_{|K}}\otimes_{\co_{|K}} L:(v_i, \alpha)\mapsto v_i\otimes \alpha,$$ which is bilinear, and the tensor product map
    $$\otimes : V_K\times L\rightarrow V_K\otimes_K L: (v, \alpha)\mapsto v\otimes \alpha.$$

    Again, by universal property, there is a unique map $$\varphi_2: V_K\otimes_K L \rightarrow \Lambda_{\co_{|K}}\otimes_{\co_{|K}} L$$ such that $\varphi_2\circ \otimes=f_2,$ which gives $\varphi_2(v_i\otimes\alpha)=v_i\otimes\alpha$. Furthermore, observe that $\varphi_1\circ\varphi_2(v_i\otimes \alpha)=v_i\otimes \alpha=\varphi_2\circ\varphi_1$, which implies $\Lambda_{\co_{|K}}\otimes_{\co_{|K}} L \cong V_K \otimes_K L$.
    
    We now claim that $V_K \otimes_K L\cong V$. To see this, consider the restriction of the natural isomorphism $f: V\otimes_L L\rightarrow V: v\otimes\alpha\mapsto \alpha v$ to
    $$f_{V_K\otimes_K L}: V_K\otimes_K L\rightarrow V: v\otimes\alpha\mapsto \alpha v.$$ Clearly, $f_{V_K\otimes_K L}$ is injective because it is a restriction of an injective function. Moreover, from universal property it follows that $f_{V_K\otimes_K L}$ is a surjective map. 

    Also, by a similar argument as above, it follows that $V\cong \Lambda\otimes_{\co} L$. In conclusion, we have the following: $$\Lambda_{\co_{|K}}\otimes_{\co_{|K}} L \cong V_K \otimes_K L\cong V\cong \Lambda\otimes_{\co} L.$$

    The above isomorphism then restricts to an isomorphism of $\co$-modules $$\Psi: \Lambda_{\co_{|K}}\otimes_{\co_{|K}} \co\rightarrow\Lambda\otimes_{\co} \co\cong\Lambda,$$ which also preserves the form $Q$. Hence $(\Lambda, Q)$ is defined over $K$.
\end{proof}

\section{Estimate on small elements}\label{s4}

Let $m\in\Z_{\geq 2}$ be even and $K$ a totally real number field. In this section, we first give an upper bound on the norm of all $\alpha\in\co_K^+$ such that $\alpha=\beta^m+\gamma$ for some $\beta,\gamma\in\co_K$. Then we show a way of expressing integers in $\co_K^+$ using $m$th powers.

\begin{lemma}\label{3,1}
    Let $K$ be a totally real number field of degree $d=[K:\Q]$ and discriminant $\Delta$. For every $\alpha\in\co_K^+$ with $\N_{K/\Q}(\alpha)>\Delta^{m/2}$ there exists $\beta\in\co_K$ such that $\alpha\succ \beta^m.$
\end{lemma}

\begin{proof}
    Let $\alpha\in\co_K^+$ with $\N_{K/\Q}(\alpha)>\Delta^{m/2}$. Then there is an $\epsilon>0$ such that $$\sigma_i(\alpha)^{1/m}-\epsilon>\Delta^{m/2}$$ for all $1\leq i\leq d$. Consider the set 
    $$B=\{(x_1,\dots,x_d)\in K_{\R} \mid |x_i|\leq \sigma_i(\alpha)^{1/m}-\epsilon\},$$ where $K_{\R}$ is the Minkowski space associated to $K$; is a compact, convex, symmetric set with the volume $$\text{Vol}(B)=\int_{-\sigma_i(\alpha)^{1/m}+\epsilon}^{\sigma_i(\alpha)^{1/m}-\epsilon}\cdots\int_{-\sigma_i(\alpha)^{1/m}+\epsilon}^{\sigma_i(\alpha)^{1/m}-\epsilon} dx_1\dots d{x_d}>2^d\Delta^{m/2}.$$ By Minkowski's theorem on lattices (see, for instance, \cite[Theorem $7.1$]{st}), there exists $0\neq \beta\in\co_K$ such that $\sigma(\beta)=(\sigma_1(\beta),\dots,\sigma_d(\beta))\in B$. We then have $\sigma_i(\beta)<\sigma_i(\alpha)^{1/m}$ for all $1\leq i\leq d$, which implies $\alpha\succ\beta^m$. 
\end{proof}

\begin{theorem}\label{exists}
    Let $K$ be a totally real number field of degree $d=[K:\Q]$ and discriminant $\Delta$. Fix a set of representatives $\mathcal{S}$ for classes of elements $\alpha\in\co_K^+$ with $\N_{K/\Q}(\alpha) \leq \Delta^{m/2}$. Then every element of $\co_K^+$ can be written as 
    \begin{equation}\label{e2}
        \sum_{\alpha\in\mathcal{S}}\alpha x_\alpha ^m+ U,
    \end{equation}
    where $U$ is a sum of $m$th powers of elements in $\co_K$.
\end{theorem}

\begin{proof} 
    Let $\gamma\in\co_K^+$ be arbitrary. If $\N_{K/\Q}(\gamma)\leq \Delta^{m/2}$, then there exists $\alpha_1\in\mathcal{S}$ such that $\gamma=\alpha_1 \beta^m$ for some $\beta\in\co_K$. Thus, by setting $x_{\alpha_1}=\beta$ and $x_\alpha=0$ for all $\alpha\neq\alpha_1$, and $U=0$, we can see that $\gamma$ can be written as \ref{e2}.

    If $\N_{K/\Q}(\gamma)>\Delta^{m/2}$, then by Lemma \ref{3,1} there exists $\beta_1\in\co_K$ such that $\gamma\succ\beta_1^m$. Write $\gamma=\beta_1^m+\delta_1$ for some $\delta_1\in\co_K^+$. If $\N_{K/\Q}(\delta_1)\leq \Delta^{m/2}$, then we can represent $\delta_1$ by $\sum_{\alpha\in\mathcal{S}}\alpha x_\alpha ^m$, and by taking $U=\beta_1^m$ we see that $\gamma$ can be written as \ref{e2}. But if $\N_{K/\Q}(\delta_1)>\Delta^{m/2}$, then we apply Lemma \ref{3,1} to $\delta_1$ again. After repeatedly applying Lemma \ref{3,1}, and considering that there are only finitely many integers in $\co_K^+$ with norms between $\Delta^{m/2}$ and $\N_{K/\Q}(\delta)$ up to multiplication by units, there exists $k\in\Z_{>0}$ such that $\gamma=\beta_1^m+\cdots+\beta_k^m+\delta_k$, where $\N_{K/\Q}(\delta_k)\leq \Delta^{m/2}$. Now, by taking $U=\beta_1^m+\cdots+\beta_k^m$ we can write $\gamma$ as \ref{e2}.
\end{proof}

\begin{rem}
     As compared to \cite[Theorem $5.4$]{kp}, Theorem \ref{exists} is not very precise but is slightly less technical.
\end{rem}

\section{Weak lifting problem for higher degree forms}\label{s5}
 Throughout this section, let $K$ be a totally real number field of degree $k=[K:\Q]$ and $\co$ an order in $K$. Recall that an order $\co\subseteq\co_K$ is a subring such that $\text{Frac }\co=K$. Consider the Minkowski embedding $\sigma: K\rightarrow \R^k$ given by $\sigma(\alpha)=(\sigma_1(\alpha),\dots,\sigma_k(\alpha))$, viewing $\sigma(K)$ as a subset of $\R^k$. Without loss of generality, we regard $\R^k$ as a topological space with the topology induced by the standard $\ell^2$ norm defined as follows: for $x=(x_1,x_2,\dots,x_k)\in\R^k$, $$\parallel x\parallel=(|x_1|^2+\cdots+|x_k|^2)^{1/2},$$ where $|.|$ is the standard absolute value on $\R$. The subset $\sigma(K)$ inherits the subspace topology. By $\mathbb{S}^1$ we mean $\mathbb{S}^1=\{x\in\R^k : \parallel x\parallel=1\}\subset\R^k.$

 Given an $m$-ic $\co$-lattice $(\Lambda, Q)$ of rank $n$. Consider the $K$-vector space $V_K=K\Lambda$. By fixing a suitable $K$-basis $v_1, v_2,\dots,v_n$ of $V_K$, every $v\in\Lambda$ can be written as $v=\gamma_1 v_1+\gamma_2 v_2+\cdots+\gamma_n v_n$ with $\gamma_i\in\co$. Hence, throughout this section, we will denote $v$ by the coordinate $v=(\gamma_1,\dots,\gamma_n)\in \Lambda$. 
\begin{proposition}\label{p1}
    Let $n\in\Z_{>0}$, $m\in\Z_{\geq 2}$ be even, $K$ a totally real number field, $\co$ an order in $K$, and $(\Lambda, Q)$ a totally positive definite $m$-ic $\co$-lattice of rank $n$. Then there exists $C\in\mathbb{Z}_{>0}$ such that for every $v=(\gamma_1,\dots,\gamma_n)\in\Lambda$ we have 

    \begin{enumerate}
        \item $Q(v)\succ \frac{1}{C}(\gamma_1^2+\gamma_2^2+\cdots+\gamma_n^2)^{\frac{m}{2}},$ and 
        \item $\house{v} ^m\leq C\house{Q(v)}.$
    \end{enumerate}
\end{proposition}

\begin{rem}
    The constant $C\in\Z_{>0}$ depends on $\Lambda,Q,K$ and the choice of $K$-basis $v_1,\dots,v_n$ for $V_K=K\Lambda$.
\end{rem}
\begin{proof}
    (1) Let $(\Lambda, Q)$ be a totally positive definite $m$-ic $\co$-lattice of rank $n$. For contradiction, suppose that for every $t\in\Z_{>0}$ there exists an embedding $\sigma : K\rightarrow \R$ and a non-zero vector $v_t=(\gamma_{t1},\dots,\gamma_{tn})\in \Lambda$ such that 
    \begin{equation*}
        \sigma\left(Q(v_t)-\frac{1}{t}(\gamma_{t1}^2+\gamma_{t2}^2+\cdots+\gamma_{tn}^2)^{\frac{m}{2}}\right)\leq 0.
    \end{equation*}
    Since there are only finitely many embeddings one $\sigma$ occurs infinitely often, let us fix that $\sigma$. Then we have 
    \begin{equation*}
        \sigma\left(Q(v_t)-\frac{1}{t}(\gamma_{t1}^2+\gamma_{t2}^2+\cdots+\gamma_{tn}^2)^{\frac{m}{2}}\right) =\sigma(Q)(\sigma(v_t))-\frac{1}{t} \left (\sigma(\gamma_{t1})^2+\cdots+\sigma(\gamma_{tn})^2\right)^{\frac{m}{2}} \leq 0.
    \end{equation*}
    Set $$w_t=(\theta_{t1},\dots,\theta_{tn})=\frac{1}{\parallel v_t\parallel}\left(\sigma(\gamma_{t1}),\dots,\sigma(\gamma_{tn})\right)=\frac{\sigma(v_t)}{\parallel v_t\parallel}.$$ Thus, we have 
    $$ \sigma(Q)\left(\parallel v_t\parallel w_t\right)-\frac{1}{t}\left(\parallel v_t\parallel ^2\theta_{t1}^2+\cdots+\parallel v_t\parallel^2\theta_{tn}^2\right)^{\frac{m}{2}}\leq 0,$$
    $$\parallel v_t\parallel^m \sigma(Q)(w_t)-\parallel v_t\parallel^m\frac{1}{t}(\theta_{t1}^2+\cdots+\theta_{tn}^2)^{\frac{m}{2}}\leq 0,$$
    $$\sigma(Q)(w_t)-\frac{1}{t}(\theta_{t1}^2+\cdots+\theta_{tn}^2)^{\frac{m}{2}}\leq 0,$$
    as $v_t$ is a non-zero.
   
    Note that $(w_t)_{t\in\Z_{>0}}$ is an infinite sequence of vectors in $\mathbb{S}^1\subset \R^n$, and since $\mathbb{S}^1$ is compact, there exists $w\in \mathbb{S}^1$ such that $w_t\rightarrow w$ as $t\rightarrow \infty$. Now, since $\sigma(Q)$ is continuous, it follows that as $t\rightarrow\infty$ we have
    $$\sigma(Q)(w)\leq 0,$$
    which implies $Q$ is not totally positive definite, a contradiction.

    (2) For every $v=(\gamma_1,\dots,\gamma_n)\in \Lambda$, we have 
    \begin{equation*}
        Q(v)\succ \frac{1}{C}(\gamma_1^2+\gamma_2^2+\cdots+\gamma_n^2)^{\frac{m}{2}}\succeq \frac{1}{C} \gamma_i^m,
    \end{equation*}
    for all $i$.
    Thus for every embedding $\sigma$ of $K$
    \begin{equation*}
        \sigma(Q(v))\geq \frac{1}{C} \sigma(\gamma_i)^m,
    \end{equation*}
    which implies 
    \begin{equation*}
        \house{Q(v)}=\max_{\sigma} |\sigma(Q(v))|\geq \frac{1}{C}
        \max_{\sigma} |\sigma(\gamma_i)^m|=\frac{1}{C}\house{\gamma_i} ^m \text{ for all }i, 
    \end{equation*}
    $$C\house{Q(v)}\geq \max_{1\leq i\leq n} \house{\gamma_i}^m=\house{v}^m$$
    as we wanted.
\end{proof}

Using Proposition \ref{p1} and following the proof of  \cite[Theorem $7$]{ky} we establish our first main result.
\begin{theorem}\label{liftingproblem}
    Let $m\in\Z_{\geq 2}$ be even, $K$ be a totally real number field, $(\Lambda, Q)$ a totally positive definite $m$-ic $\co_K$-lattice, and $\alpha\in\co_K^+.$ There exist at most finitely many totally real number fields $L/K$ of degree $d=[L : \Q]$ such that $(\Lambda\otimes_{\co_K}\co_L, Q_L)$ represents all $\alpha\co_L^+.$
\end{theorem}

\begin{proof}
    Fix an integral basis $\omega_1, \omega_2,\dots, \omega_k$ for $K,$ and $\alpha\in\co_K^+.$ Let $L/K$ be an extension of totally real number field of degree $d=[L:\Q].$  Choose a real number $X>1$ such that $$X>\max \left(\max_{1\leq i\leq k}\house{\omega_i}, 3C\N_{K/\Q}(\alpha)\right),$$ where $C$ is the constant from Proposition \ref{p1}.
    
     Now, consider the proper subfield $F\subset L$ generated by $\beta\in\co_L$ with $\house{\beta}<X.$ Since $X$ was chosen such that $\house{\omega_i}<X,$ we have $\omega_i\in F$ for all $1\leq i\leq k,$ hence $K\subset F.$ By Northcott property, there are only finitely many choices for $F.$
     
     \medskip
     Also, by Northcott property there exists $X\leq X_1=\min \left(\house{\beta} \mid \beta\in \co_L\setminus \co_F \right)$. Let $\beta_1\in \co_L\setminus \co_F$ be such that $\house{\beta_1}=X_1,$ choose an smallest integer $t$ such that $\beta_1+t\succ 0.$ Then we have 
     $$\house{\beta_1 + t}\leq \house{t}+\house{\beta_1}<2\house{\beta_1}+1\leq 3\house{\beta_1}.$$

     Suppose $(\Lambda\otimes_{\co_K} \co_L, Q_L)$ represents all $\alpha\co_L^+,$ in particular there exists $v\in \Lambda\otimes_{\co_K}\co_L$ such that $Q_L(v)=\alpha(\beta_1+t).$ If $v=(v_1,\dots,v_n)\in F^n$, since $(\Lambda,Q)$ is an $m$-ic $\co_K$-lattice and $K\subset F$, we have $Q_L(v)=\alpha(\beta_1+t)\in F$, which implies $\beta_1\in F,$ a contradiction. Thus, there exists $1\leq i\leq n$ such that $v_i\notin F$, therefore we have $\house{v}\geq \house{v_i}\geq X_1.$

     By Proposition \ref{p1} we further have that 
     $$X_1^m\leq \house{v_i}^m\leq C \house{Q(v)}=C \house{\alpha (\beta_1+t)}\leq C \house{\alpha}\house{\beta_1+t}<3C \N_{K/\Q}(\alpha)\house{\beta_1}=3CX_1\N_{K/\Q}(\alpha),$$
     which implies
     $$X\leq X_1^{m-1}<3C\N_{K/\Q}(\alpha).$$ This contradicts our choice of $X$.
\end{proof}

\begin{theorem}\label{fail.ST}
    Let $n,d\in\Z_{>0}$, $m\in\Z_{\geq 2}$ be even, $K$ a totally real number field, $(\Lambda, Q)$ a totally positive definite $m$-ic $\co_K$-lattice of rank $n$, and $\alpha\in\co_K$. There are only finitely many $L/K$ of degree $d=[L:K]$ such that $\alpha$ is represented by $(\Lambda, Q)\otimes_{\co_K} \co_L$ but is not represented by $(\Lambda, Q)\otimes_{\co_K} \co_{L'}$ for every proper subfield $L'$ of $L$ containing $K$.
\end{theorem}

\begin{proof}
    Recall that if $v_1,v_2,\dots,v_n$ is a basis for $V_K=K\Lambda$, then every $v\in\Lambda$ can be written as $v=\gamma_1v_1+\cdots+\gamma_nv_n$, where $\gamma_i\in\co_K$. 

    Let $L/K$ be a totally real extension of $K$ of degree $d=[L:K]$. Observe that $\Lambda\otimes_{\co_K}\co_L\subseteq \co_L v_1+\cdots+\co_L v_n$. 

    Suppose that $\alpha\in\co_K$ is represented by $(\Lambda, Q)\otimes_{\co_K} \co_L$. Then there exists $v\in\Lambda\otimes_{\co_K} \co_L$ such that $Q_L(v)=Q_L(\gamma_1v_1+\cdots+\gamma_nv_n)=\alpha$, where $\gamma_i\in\co_L$ for all $1\leq i\leq n$. Thus, we have $K(\gamma_1,\gamma_2,\dots,\gamma_n)\subseteq L$. Further, if $\alpha$ is not represented by $(\Lambda, Q)\otimes_{\co_K} \co_{L'}$ for every proper subfield $L'$ of $L$ containing $K$, then $L=K(\gamma_1,\gamma_2,\dots,\gamma_n)$. By Proposition \ref{p1}, it then follows that $$\house{v}^m<C\house{\alpha}.$$ Consequently, by Northcott property it follows that there are only finitely many choice for $\gamma_i s$, and hence only finitely many such field $L$ may exist. 
\end{proof}

\section{(Non-)universality of $m$-ic lattices over infinite extensions}\label{s6}

We begin by extending Proposition \ref{p1} from totally real number fields to totally real infinite extensions of $\Q$.

\begin{proposition}\label{p2}
    Let $n\in\Z_{>0}$, $m\in\Z_{\geq 2}$ be even, $K$ be totally real extension of $\Q$, $\co$ an order in $K$, and $(\Lambda, Q)$ be totally positive definite $m$-ic $\co$-lattice of rank $n$. Then there is a constant $\kappa=\kappa(Q)\in \Z_{>0}$ such that for every $v\in \Lambda$ we have
    $$\house{v} ^m\leq \kappa\house{Q(v)}.$$
\end{proposition}

\begin{proof}
Let $(\Lambda, Q)$ be an $m$-ic $\co$-lattice of rank $n$. By Proposition \ref{definedover}, $(\Lambda, Q)$ is defined over some number field $F\subset K$, i.e., there exists an $m$-ic $\co_{|F}$-lattice $(\Lambda_{\co_{|F}}, Q_F)$ such that $(\Lambda, Q)$ is a scalar extension of $(\Lambda_{\co_{|F}}, Q_F)$.

Define $V_F=F\Lambda_{\co_{|F}}$ is a $F$-vector space of dimension $n$. Let $v_1',v_2',\dots,v_n'$ be an $F$-basis for $V_F$ and $u_1,\dots,u_k$ be a set of generators of $\Lambda_{\co_{|F}}$ (as $\co_F$-submodule). Thus, we can write each $u_i=\lambda_{i1}v_1'+\cdots+\lambda_{in}v_n',$ where $\lambda_{ij}\in F$. Since $\text{Frac }\co_{|F}=F$, we have $\lambda_{ij}=a_{ij}/b_{ij}$ where $a_{ij}, b_{ij}\in\co_{|F}$. In fact, by rescaling, we may choose a uniform denominator, i.e., $\lambda_{ij}=a_{ij}/b$ with $a_{ij}, b\in\co_{|F}$. Note that if $v_1,v_2,\dots,v_n$ is an $F$-basis for $V_F$, then $v_1=(1/b)v_1', v_2=(1/b)v_2',\dots,v_n=(1/b)v_n'$ also form an $F$-basis for $V_F$. Thus, we may assume $\lambda_{ij}\in\co_{|F}$.

Let $v\in\Lambda$ be an arbitrary vector. Since $\Lambda$ is an scalar extension of $\Lambda_{\co_{|F}}$ we can write $v$ as $v=\mu_1w_1+\mu_2w_2+\cdots+\mu_lw_l$, where $\mu_i\in \co$ and $w_i\in\Lambda_{\co_{|F}}$; and $w_i=\gamma_{i1}u_1+\cdots+\gamma_{ik}u_K$, where $\gamma_{ij}\in\co_{|F}$. In total, we can write $v=\gamma_1v_1+\gamma_2v_2+\cdots+\gamma_nv_n$ with $\gamma_i\in \co$. If $Q(v)=\alpha\in\co^+$, then $\gamma_1,\gamma_2,\dots,\gamma_n,\alpha$ lies in the number field $F\subseteq F'=F(\gamma_1,\gamma_2,\dots,\gamma_n,\alpha)\subset K$. Then, by Proposition \ref{p1} there is a constant $\kappa\in\Z_{>0}$ such that $$\house{v}^m\leq \kappa \house{Q(v)}=\kappa\house{\beta}.$$
This completes the proof.
\end{proof}

\begin{corollary}
    Let $K$ be totally real extension of $\Q$, $\alpha\in\co^+$, and $(\Lambda, Q)$ be totally positive definite $m$-ic $\co$-lattice. If $\co^+$ has the Northcott property, then we have 
    \begin{equation*}
        \#\{v\in\Lambda \mid Q(v)=\alpha\}<\infty.
    \end{equation*}
    
\end{corollary}

\begin{proof}
    Suppose $(\Lambda, Q)$ has rank $n$. By Proposition \ref{p2}, for every $v\in \Lambda$ $$\house{v} ^m\leq \kappa\house{Q(v)}.$$ But then, by Northcott property we have 
    $$\#\{v\in\Lambda \mid Q(v)=\alpha\}\leq \#\{v\in\Lambda \mid \house{v}\leq \kappa \house{\alpha}\}<\infty.\qedhere$$

\end{proof}

The above corollary justifies our definition of positive definiteness. For example, consider the form $Q(x,y)=(x^2-2y^2)^m$ over the field $\Q$, then $Q(x,y)=1$ has infinitely many solution $(x,y)\in\Z^2$ as the Pell's equation $x^2-2y^2=1$ has infinitely many solution in integers $x,y.$

\begin{theorem}\label{thm:infinite}
    Let $m\in\Z_{\geq 2}$ be even, $K$ be a totally real infinite extension of $\Q$, $\delta\in\co_K^+$, and $\co$ be an order in $K$. If there exists an $m$-ic $\co$-lattice that represents all $\delta\co^+$, then $\co^+$ does not have the Northcott property. 
\end{theorem}

\begin{proof}
    Let $\delta\in\co_K^+$ be fixed. Assume that there exists an $m$-ic $\co$-lattice $(\Lambda, Q)$ that represents all $\delta\co^+$. By Proposition \ref{definedover}, there is a number field $F\subset K$ and an $m$-ic $\co_{|F}$-lattice $(\Lambda_{\co_{|F}}, Q_F)$ such that $(\Lambda, Q)$ is a scalar extension of $(\Lambda_{\co_{|F}}, Q_F)$ to $\co$. Let $\alpha\in\co^+$ be any element, then there is a vector $v\in\Lambda$ such that $Q(v)=\delta\alpha$. Note that, as shown in the proof of Proposition \ref{p2}, with a suitable $F$-basis $v_1,\dots,v_n$ of $V_F=F\Lambda_{\co_{|F}}$ we may write $v=\alpha_1v_1+\alpha_2v_2+\cdots+\alpha_nv_n$, where $\alpha_1,\alpha_2,\dots,\alpha_n \in \co$ and depends on the choice of $\alpha$.

    Define $$\beta_i=\alpha_i+\lfloor \house{\alpha_i}\rfloor +1$$ for all $1\leq i\leq n$. Observe that $\beta_i\in\co^+$ since $\lfloor\house{\alpha_i}\rfloor\succ -(1+\alpha_i)$ holds, thus $\delta\beta_i\in\delta\co^+$.

    Now, define $A_0=\{\delta\alpha\}$ and $A_1=\{\delta\beta_1,\delta\beta_2,\dots,\delta\beta_n\}.$ Moreover, from Proposition \ref{p2} it follows that 
    \begin{equation*}
        \house{\delta\beta_i}=\house{\delta}\house{\alpha_i+\lfloor \house{\alpha_i}\rfloor +1}<\house{\delta}\left(2\house{\alpha_i}+1<2\kappa^{1/m}\house{\alpha}^{1/m}+1\right)<\house{\delta}(2\kappa^{1/m}+1)\house{\alpha}^{1/m}.
    \end{equation*}
    Recursively, we construct the sets $A_j\subseteq\co^+$ for all $j\geq 2$ as follows: for $\delta\beta\in A_{j-1}$ there is $v=(\alpha_1,\alpha_2,\dots,\alpha_n)\in\Lambda$ such that $Q(v)=\delta\beta$. We then define $\beta_i=\alpha_i+\lfloor \house{\alpha_i}\rfloor +1$ for all $i$. Clearly, as discussed earlier, we have $\delta\beta_i\in\delta\co^+$, and $$\house{\delta\beta_i}=\house{\delta}\house{\alpha_i+\lfloor \house{\alpha_i}\rfloor +1}<\house{\delta}\left(2\house{\alpha_i}+1<2\kappa^{1/m}\house{\beta}^{1/m}+1\right)<\house{\delta}(2\kappa^{1/m}+1)\house{\beta}^{1/m}.$$ Now, take $A_j=\{\delta\beta_1,\dots,\delta\beta_n \mid \delta\beta \in A_{j-1}\}\subseteq\delta\co^+$. 
    \medskip
    
    As $j\rightarrow\infty$, $\house{\delta\beta}<\house{\delta}(2\kappa^{1/m}+1)^{m/m-1}$ for all $\beta\in A_j$. Therefore, for sufficiently large $j$ we have $$\house{\delta\beta}<\house{\delta}(2\kappa^{1/m}+1)^{m/m-1}\text{ for all }\delta\beta\in A_j.$$

    For contradiction, suppose that $\co^+$ has the Northcott property. Then there exist only finitely many $\gamma\in\co^+$ such that $\house{\gamma}<\house{\delta}(2\kappa^{1/m}+1)^{m/m-1}$. Let $F'=F(\gamma\in\co^+ \mid \house{\gamma}<\house{\delta}(2\kappa^{1/m}+1)^{m/m-1})\subset K$ be the number field generated by $\gamma\in\co^+$ with $\house{\gamma}<\house{\delta}(2\kappa^{1/m}+1)^{m/m-1}$. Given $\alpha\in\co^+$, for sufficiently large $j$ we have $A_j\subseteq \co_{F'}^+$. Consequently, by construction, $A_0\subseteq\co_{F'}^+$, i.e., $\delta\alpha\in\co_{F'}^+$. Thus $\delta\co^+\subseteq\co_{F'}^+$ which implies $K\subseteq F'$, a contradiction to the fact that $K$ is an infinite extension of $\Q$. Therefore, $\co^+$ does not possess the Northcott property.
\end{proof}

For a number field $K$ and a positive integer $d$, let $K^{(d)}$ be the compositum of all Galois extensions of $K$ of degree at most $d$, and $K_{\text{ab}}^{(d)}$ be the maximal abelian extension $F/K$ with $K\subseteq F \subseteq K^{(d)}$. Bombieri and Zannier \cite[Theorem $1$]{bz} showed that $K_{\text{ab}}^{(d)}$ has the Northcott property with respect to Weil height for any $d\in\Z_{>0}$. Hence, it follows that $K_{\text{ab}}^{(d)}$ has the Northcott property with respect to house. 

\begin{corollary}\label{c2}
    Let $m\in\Z_{\geq 2}$ be even, $d\in\Z_{>0}$ be an integer, $K$ be a number field, and $K\subseteq L\subseteq K_{\text{ab}}^{(d)}$ be a totally real field of infinite degree $[L:\Q]$. Then there is no universal $m$-ic $\co$-lattice, where $\co$ is an order in $L$.
\end{corollary}

\begin{proof}
    Since $K_{\text{ab}}^{(d)}$ has the Northcott property, it follows that any subfield $L\subseteq K_{\text{ab}}^{(d)}$ also possesses the Northcott property. Consequently, $\co^+$ also has the Northcott property. Now, the statement follows immediately from Theorem \ref{thm:infinite}. 
\end{proof}

A positive integer $n\in\Z_{>0}$ is said to be an \emph{abelian number (or abelianness-forcing number)} if every group of order $n$ is abelian. From \cite[pp. $187$, Exercise $24$]{dummit} it follows that $n\in\Z_{>0}$ is an abelian number if and only if $n$ is cube-free and there is no prime power $p^k|n$ with $k\geq 1$, such that $p^k\equiv 1 \pmod q$ for some prime $q|n$. For example, if $p$ is a prime number then $p$ and $p^2$ are abelian numbers; further, if $p$ is odd and $q>p$ is a prime such that $q\not\equiv 1 \pmod p$, then $pq$ is an abelian number.

Now, using Corollary \ref{c2}, we provide a large class of examples of totally real fields that have no universal $m$-ic forms. Let us denote by $\Q^{[d]}$ the composititum of all totally real Galois extensions of $\Q$ of degree exactly $d$.
\begin{corollary}
    Let $n\in\Z_{>0}$ be an abelian number. If $K$ is an infinite extension of $\Q$ contained in $\Q^{[n]}$, and $\co$ is an order in $K$, then there are no universal $m$-ic $\co$-lattices.
\end{corollary}

\begin{proof}
    Since $n$ is an abelian number, it follows that $\Q^{[n]}$ is abelian, hence a subfield of $\Q_{\text{ab}}^{(n)}$. If $K$ is an infinite extension of $\Q$ contained in $\Q^{[n]}$, then $K\subseteq \Q_{\text{ab}}^{(n)}$. Now, from Corollary \ref{c2}, the assertion follows.
\end{proof}

\bibliographystyle{amsalpha}
	
\bibliography{Citation}

\providecommand{\bysame}{\leavevmode\hbox to3em{\hrulefill}\thinspace}
\providecommand{\MR}{\relax\ifhmode\unskip\space\fi MR }
\providecommand{\MRhref}[2]{%
  \href{http://www.ams.org/mathscinet-getitem?mr=#1}{#2}
}
\providecommand{\href}[2]{#2}
\begin{thebibliography}{HLYZ21b}

\bibitem[BH11]{bh}
M.~Bhargava and J.~Hanke, \emph{Universal quadratic forms and the 290-theorem}, In preprint, 2011.

\bibitem[BK15]{bk1}
V.~Blomer and V.~Kala, \emph{Number fields without {$n$}-ary universal quadratic forms}, Math. Proc. Cambridge Philos. Soc. \textbf{159} (2015), no.~2, 239--252. \MR{3395370}

\bibitem[BK18]{bk2}
\bysame, \emph{On the rank of universal quadratic forms over real quadratic fields}, Doc. Math. \textbf{23} (2018), 15--34. \MR{3846065}

\bibitem[BZ01]{bz}
E.~Bombieri and U.~Zannier, \emph{A note on heights in certain infinite extensions of {$\Bbb Q$}}, Atti Accad. Naz. Lincei Cl. Sci. Fis. Mat. Natur. Rend. Lincei (9) Mat. Appl. \textbf{12} (2001), 5--14 (2002). \MR{1898444}

\bibitem[DF04]{dummit}
D.~S. Dummit and R.~M. Foote, \emph{Abstract algebra}, third ed., John Wiley \& Sons, Inc., Hoboken, NJ, 2004. \MR{2286236}

\bibitem[DKKY24]{dkky}
N.~Daans, V.~Kala, J.~Kr\'{a}sensk\'{y}, and P.~Yatsyna, \emph{Failures of integral springer's theorem}, In preprint, https://arxiv.org/pdf/2404.12844, 2024.

\bibitem[DKM23]{dkm}
N.~Daans, V.~Kala, and S.~H. Man, \emph{Universal quadratic forms and northcott property of infinite number fields}, In preprint, https://arxiv.org/pdf/2308.16721, 2023.

\bibitem[FP09]{ab5}
R.~W. Fitzgerald and S.~Pumpl\"{u}n, \emph{Norm principles for forms of higher degree permitting composition}, Comm. Algebra \textbf{37} (2009), no.~11, 3851--3860. \MR{2573223}

\bibitem[Har75]{harrison}
D.~K. Harrison, \emph{A {G}rothendieck ring of higher degree forms}, J. Algebra \textbf{35} (1975), 123--138. \MR{379370}

\bibitem[HLYZ21a]{ab1}
H-L. Huang, H.~Lu, Y.~Ye, and C.~Zhang, \emph{Diagonalizable higher degree forms and symmetric tensors}, Linear Algebra Appl. \textbf{613} (2021), 151--169. \MR{4195104}

\bibitem[HLYZ21b]{huang}
\bysame, \emph{Diagonalizable higher degree forms and symmetric tensors}, Linear Algebra Appl. \textbf{613} (2021), 151--169. \MR{4195104}

\bibitem[JP10]{pumpl}
A.~Johnson and S.~Pumpl\"{u}n, \emph{Anisotropic forms of higher degree under finite dimensional field extensions}, Comm. Algebra \textbf{38} (2010), no.~6, 2064--2078. \MR{2675522}

\bibitem[Kal16a]{k}
V.~Kala, \emph{Norms of indecomposable integers in real quadratic fields}, J. Number Theory \textbf{166} (2016), 193--207. \MR{3486273}

\bibitem[Kal16b]{k2}
\bysame, \emph{Universal quadratic forms and elements of small norm in real quadratic fields}, Bull. Aust. Math. Soc. \textbf{94} (2016), no.~1, 7--14. \MR{3539315}

\bibitem[Kal23]{survey}
V.~Kala, \emph{Universal quadratic forms and indecomposables in number fields: a survey}, Commun. Math. \textbf{31} (2023), no.~2, 81--114. \MR{4621253}

\bibitem[KP24]{kp}
V.~Kala and O.~Prakash, \emph{There is no $290$-theorem for higher degree forms}, In preprint, https://arxiv.org/pdf/2405.10660, 2024.

\bibitem[KT23]{kt}
V.~Kala and M.~Tinkov\'{a}, \emph{Universal quadratic forms, small norms, and traces in families of number fields}, Int. Math. Res. Not. IMRN (2023), no.~9, 7541--7577. \MR{4584708}

\bibitem[KW90]{kw}
J.~M. Kubina and M.~C. Wunderlich, \emph{Extending {W}aring's conjecture to {$471,600,000$}}, Math. Comp. \textbf{55} (1990), no.~192, 815--820. \MR{1035936}

\bibitem[KY23]{ky}
V.~Kala and P.~Yatsyna, \emph{On {K}itaoka's conjecture and lifting problem for universal quadratic forms}, Bull. Lond. Math. Soc. \textbf{55} (2023), no.~2, 854--864. \MR{4581328}

\bibitem[Mor07]{ab2}
P.~J. Morandi, \emph{A {S}pringer theorem for higher degree forms}, Math. Z. \textbf{256} (2007), no.~1, 221--228. \MR{2282266}

\bibitem[Nor49]{north}
D.~G. Northcott, \emph{An inequality in the theory of arithmetic on algebraic varieties}, Proc. Cambridge Philos. Soc. \textbf{45} (1949), 502--509. \MR{33094}

\bibitem[O'M63]{omeara}
O.~T. O'Meara, \emph{Introduction to quadratic forms}, Die Grundlehren der mathematischen Wissenschaften, Band 117, Springer-Verlag, Berlin-G\"{o}ttingen-Heidelberg; Academic Press, Inc., Publishers, New York, 1963. \MR{152507}

\bibitem[Pum09]{ab4}
S.~Pumpl\"{u}n, \emph{{$u$}-invariants for forms of higher degree}, Expo. Math. \textbf{27} (2009), no.~1, 37--53. \MR{2503042}

\bibitem[Pum13]{ab6}
\bysame, \emph{Round forms of higher degree}, Results Math. \textbf{63} (2013), no.~1-2, 657--674. \MR{3009712}

\bibitem[Pum20]{ab3}
\bysame, \emph{Diagonal forms of higher degree over function fields of {$p$}-adic curves}, Int. J. Number Theory \textbf{16} (2020), no.~1, 161--172. \MR{4059129}

\bibitem[Sie21]{si1}
C.L. Siegel, \emph{Darstellung total positiver {Z}ahlen durch {Q}uadrate}, Math. Z. \textbf{11} (1921), no.~3-4, 246--275. \MR{1544496}

\bibitem[Sie23]{si2}
\bysame, \emph{Additive {T}heorie der {Z}ahlk\"{o}rper. {II}}, Math. Ann. \textbf{88} (1923), no.~3-4, 184--210. \MR{1512127}

\bibitem[Sie44]{siegel}
\bysame, \emph{Generalization of {W}aring's problem to algebraic number fields}, Amer. J. Math. \textbf{66} (1944), 122--136. \MR{9778}

\bibitem[ST02]{st}
I.~Stewart and D.~Tall, \emph{Algebraic number theory and {F}ermat's last theorem}, third ed., A K Peters, Ltd., Natick, MA, 2002. \MR{1876804}

\bibitem[VV16]{vv16}
X.~Vidaux and C.~R. Videla, \emph{A note on the {N}orthcott property and undecidability}, Bull. Lond. Math. Soc. \textbf{48} (2016), no.~1, 58--62. \MR{3455748}

\bibitem[VW02]{vw}
R.~C. Vaughan and T.~D. Wooley, \emph{Waring's problem: a survey}, Number theory for the millennium, {III} ({U}rbana, {IL}, 2000), A K Peters, Natick, MA, 2002, pp.~301--340. \MR{1956283}

\bibitem[Yat19]{y}
Pavlo Yatsyna, \emph{A lower bound for the rank of a universal quadratic form with integer coefficients in a totally real number field}, Comment. Math. Helv. \textbf{94} (2019), no.~2, 221--239. \MR{3942784}

\end{thebibliography}

\end{document}